\documentclass [12pt,a4paper,reqno]{amsart}
\usepackage{amsmath,amsthm,amscd,latexsym}
\usepackage{amsfonts}
\allowdisplaybreaks

\chardef\bslash=`\\ 

\theoremstyle{definition}

\newtheorem*{remark*}{Remarks}
\newtheorem*{defn*}{Definition}

\theoremstyle{remark}

\newcommand{\wt}{\widetilde}
\newcommand{\wh}{\widehat}
 \renewcommand{\sectionmark}[1]{}

 \usepackage{amsfonts}
\newcommand{\field}[1]{\mathbb{#1}}
\newcommand{\C}{\field{C}}
\newcommand{\R}{\field{R}}
\newcommand{\D}{\field{D}}
\newcommand{\T}{\mathbf{T}}
\newcommand{\hC}{\widehat{\field{C}}} 
\newcommand{\B}{\mathbf{B}}

\newcommand{\Belt}{\operatorname{Belt}}
\newcommand{\dist}{\operatorname{dist}}

\newcommand{\x} {\mathbf x}

\begin{document}

\title{Complex rigidity of Teichm\"{u}ller spaces}

\author{Samuel L. Krushkal}

\begin{abstract} We outline old and new results concerning 
the well-known problems in the Teichm\"{u}ller space theory, i.e.,   
whether these spaces are starlike in the Bers holomorphic embedding 
and whether any Teichm\"{u}ller space of dimension greater than 
$1$ is biholomorhically equivalent to bounded convex domain in a  
complex Banach space. 
\end{abstract}

\date{\today\hskip4mm ({\tt ComplexRigidT.tex})}

\maketitle

\bigskip


\medskip


\bigskip

\markboth{Samuel L. Krushkal}{Complex rigidity of Teichm\"uller spaces}
\pagestyle{headings}

\bigskip
\centerline{\bf 1. INTRODUCTORY REMARKS}  

It is well-known that the Teichm\"{u}ller spaces with their
canonical complex structure are pseudo-convex. Moreover, all finite
dimensional Teichm\"{u}ller spaces are Runge domains, hence
polynomially convex. 

The folowing two longstanding problems relate to geometric 
convexity of these spaces. 
 
\bigskip\noindent
{\bf 1}. {\it For an arbitrary finitely or infinitely generated 
Fuchsian
group $\Gamma$, is the Bers embedding of its Teichm\"{u}ller 
space $\mathbf T(\Gamma)$ starlike? }
  
\bigskip\noindent
{\bf 2}. {\it Is any finite or infinite  dimensional Teichm\"uller 
space 
of dimension greater than $1$ biholomorphically equivalent 
to bounded convex domain in a complex \linebreak Banach space $X$ (of 
the same dimension as $\mathbf T(\Gamma)$)?} 

\bigskip
The first problem was stated in among other open problems 
on Teichm\"uller spaces and Kleinian groups in  
the book \cite{BK} of 1974, collected by Abikoff. 

The second problem  was posed for the finite dimensional spaces 
by Royden and for the
universal Teichm\"{u}ller space by Sullivan. It relates to Tukia's
result \cite{Tu} which explicitly yields a real analytic
homeomorphism of the universal Teichm\"{u}ller space $\T =
\T(\mathbf 1)$ onto a convex domain in a real Banach space. 

The aim of this paper is to outline old and recent results 
obtained in solving these problems.

\bigskip\bigskip 
\centerline{\bf 2. TEICHM\"{U}LLER SPACES ARE NOT STARLIKE} 

\bigskip\noindent
{\bf 1}. \ First recall that the Bers embedding represents 
the space $\mathbf T(\Gamma)$ 
as a bounded domain formed by the Schwarzian derivatives
$$
S_w = \Bigl(\frac{w^{\prime\prime}}{w^\prime}\Bigr)^\prime 
- \frac{1}{2} \Bigl(\frac{w^{\prime\prime}}{w^\prime}\Bigr)^2
$$
of holomorphic univalent functions $w(z)$ in the lower half-plane 
$U^* = \{z: \Im z < 0\}$) (or in the disk) admitting quasiconformal 
extensions to
the Riemann sphere $\widehat{\C} = \C \cup \{\infty\}$ compatible 
with the group $\Gamma$ acting on $U^*$.

It was shown in \cite{Kr1} that the universal Teichm\"{u}ller space 
$\mathbf T = \mathbf T(\mathbf 1)$ has
points which cannot be joined to a distinguished point even by
curves of a considerably general form, in particular, by polygonal
lines with the same finite number of rectilinear segments. The proof
relies on the existence of conformally rigid domains established by
Thurston in \cite{Th} (see also \cite{As}).

This implies, in particular, that the universal Teichm\"{u}ller space 
is not starlike with
respect to any of its points, and there exist points 
$\varphi \in \mathbf T$
for which the line interval $\{t \varphi: 0 < t < 1\}$ contains the
points from $\mathbf B \setminus \mathbf S$, where 
$\mathbf B = \mathbf B(U^*)$ is the
Banach space of hyperbolically bounded holomorphic functions in the
half-plane $U^*$ with norm
\begin{equation}\label{1}
\|\varphi\|_{\mathbf B} = 4 \sup_{U^*} y^2 |\varphi (z)|
\end{equation} 
and $\mathbf S$ denotes the set of all Schwarzian derivatives of
univalent functions on $U^*$. All $\varphi$ with finite norm (1) 
determine holomorphic functions on $U^*$ 
(as solutions of the Schwarz differential equation $S_w = \varphi$)  
which are only locally univalent.

Toki \cite{To} extended the result on the nonstarlikeness of the
space $\mathbf T$ to \linebreak Teichm\"{u}ller spaces of Riemann 
surfaces that contain
hyperbolic disks of arbitrary large radius, in particular, for the
spaces corresponding to Fuchsian groups of second kind. The crucial
point in the proof of \cite{To} is the same as in \cite{Kr1}

On the other hand, it was established in \cite{Kr3} that also all
finite dimensional Teichm\"{u}ller spaces $\mathbf T(\Gamma)$ 
of high enough dimensions
are not starlike. It seems likely that this property must hold for
all Teichm\"{u}ller spaces of dimension at least two.

The non-starlikeness causes obstructions to some problems in the 
Teichm\"{u}ller space theory and its applications to geometric 
complex analysis.

\bigskip\noindent 
{\bf 2}. \ There is also a simpler proof that the universal 
Teichm\"{u}ller 
space is not starlike. This proof, given recently in \cite{Kr5},  
provides  explicitly the functions which violate this property. 
Its underlying geometric features are  completely different 
and involve the Abikoff-Bers-Zhuravlev theorem which yields that 
the domain $\mathbf T$ has a 
common boundary with its complementary domain in the space 
$\mathbf B$ (see \cite{Ab}, \cite{Be}, \cite{Zh}). 

It is technically more convenient to deal here with univalent functions 
in the upper half-plane $U = \{z = x + iy: \ y > 0\}$ denoting by 
$\mathbf B$ the corresponding space $\mathbf B(U)$ of hyperbolically 
bounded holomorphic functions in $U$.

Let $P_n$ be a convex rectilinear polygon with the finite vertices
$A_1, A_2, ... \ , A_n$, and let the interior angle at the vertex
$A_j$ be equal to $\pi \alpha_j$; then $0 < \alpha_j < 1$ and
$\sum\limits_{j=1}^n \alpha_j = n -2$.

The conformal map of $U$ onto $P_n$ is represented by the
Schwarz-Christoffel integral
\begin{equation}\label{2}
f_{*}(z) = d_1 \int\limits_0^z (\xi - a_1)^{\alpha_1 - 1} (\xi -
a_2)^{\alpha_2 - 1} ... (\xi - a_n)^{\alpha_n - 1} d \xi + d_0,
\end{equation}
where $a_j = f_{*}^{-1}(A_j) \in \R, \ a_1 < a_2 < \dots < a_n$, and
$d_0, d_1$ are the complex constants. Its logarithmic derivative
$b_f= f^{\prime\prime}/f^\prime$ is of the form 
$$ 
b_{f_{*}}(z) = \sum\limits_1^n (\alpha_j - 1)/(z - a_j),   
$$ 
and its Schwarzian
$$
S_{f_{*}}(z) = \sum\limits_1^n \frac{C_j}{(z - a_j)^2} -
\sum\limits_{j,l=1}^n \frac{C_{jl}}{(z - a_j)(z - a_l)},
$$
where 
$$
C_j = \alpha_j - 1 - \frac{1}{2} (\alpha_j - 1)^2 < 0, 
\quad C_{jl} = (\alpha_j - 1)(\alpha_l - 1) > 0.
$$
We normalize $f_{*}$ letting $a_1 =0, \ a_2 = 1$ and fixing $a_n <
\infty$; then $f_{*}(\infty)$ is an inner point of the edge $A_n A_1$.

Since the boundary $\partial P_n$ is a quasicircle, the function (2)
admits a quasiconformal extension onto the lower half-plane $U^*$, hence
$S_{f_{*}} \in \mathbf T$.

Denote by $r_0$ the positive root of the equation
 \begin{equation}\label{3}
\frac{1}{2} \Bigl[\sum\limits_1^n (\alpha_j - 1)^2 
+ \sum\limits_{j,l=1}^n
(\alpha_j - 1)(\alpha_l - 1) \Bigr] r^2 - \sum\limits_1^n 
(\alpha_j - 1) \ r - 2 = 0,
\end{equation}
and let
$$
S_{f_{*},t} = t b_{f_{*}} - \frac{t^2}{2} b_{f_{*}}^2, \quad t > 0.
$$

\bigskip\noindent 
{\bf Theorem 1}. {\em For any convex polygon $P_n$, the
Schwarzians $r S_{f_{*},r_0}$ and $S_{f_{*},r}$ with $0 <r < r_0$
define the univalent on $U$ functions, and the corresponding 
harmonic Beltrami coefficients 
$\nu_r(z) = - (r/2) y^2 S_{f_{*},r_0}(\overline z)$ 
and $\nu_r(z) = -(1/2) y^2 S_{f_{*},r}(\overline z)$
of their quasiconformal extensions to the lower half-plane $U^*$ are 
extremal
(have minimal $L_\infty$-norm). Hence, for some $r$ between $r_0$ and
$1$, the Schwarzians $r S_{f_{*},r_0}$ and $S_{f_{*},r}$ are the
outer points of $\mathbf T$.}

\bigskip 
Note that for $r < r_0$ the solutions $w_r$ of each equation
 \begin{equation}\label{4}
w^{\prime\prime}/w^\prime)^\prime - (w^{\prime\prime}/w^\prime)^2/2
= \varphi_r(z), \quad z \in U,
\end{equation}
with $\varphi_r = r S_{f_{*},r_0}$ and $\varphi_r = S_{f_{*},r}$ 
map $U$ conformally onto the quasidisks (either bounded or not), which can
be regarded as the analytic polygons with vertices 
$w_r(a_1), ... \ , w_r(a_n)$, 
whose boundary consists either of $n$ real analytic
arcs with nonzero intersection angles or else of arcs of spirals,
which are analytic in their interior points.

This theorem yields, in particular, that any such $w_r$ does not admit
extremal quasiconformal extensions of Teichm\"{u}ller type.

The coefficients $\nu_r$ define the Ahlfors-Weill quasiconformal 
extension 
of $w_r$ to the lower half-plane $U^*$, and 
$$ 
\|\nu_r\|_\infty = \frac{1}{2} \|\varphi_r\|_{\mathbf B} < 1   
$$ 
(provided that $\|\varphi_r\|_{\mathbf B} < 2$). 

The proof of Theorem 1 reveals an interesting  
connection between harmonic Beltrami coefficients and the Grunsky 
coefficient inequalities (first established in \cite{Kr4}).

\bigskip\noindent
{\bf 3}. \ 
Note that non-starlikeness of the universal Teichm\"{u}ller space is 
in fact the main step in the
proof of most of the results mentioned in the beginning of this paper. 
By appropriate
approximation, this property was extended to the spaces 
$\mathbf T(\Gamma)$ of
sufficiently large dimensions. So Theorem 1  has the same
corollaries. For example, the arguments in its proof provide 
simultaneously
non-starlikeness of the space $\mathbf T$ in Becker's holomorphic 
embedding which represents 
this space as a bounded domain in the Banach space of holomorphic  
functions $\psi$ in the disk $\Delta^* = \{|z| > 1\}$ with norm
$\|\psi\|=\sup_{\Delta^*} (|z|^2-1) z \psi(z)|$. The points of 
this domain are the logarithmic derivatives 
$\psi_f = f^{\prime\prime}/f^\prime$ of univalent functions 
$f(z) = z + b_0 + b_1 z^{-1} + \dots$ in $\Delta^*$.

\bigskip\noindent
{\bf 4}. \ 
As an {\bf example}, consider the rectangles $P_4$. For any 
rectangle, 
all $\alpha_j = 1/2$, hence the equation (4) assumes the form
$$
\frac{5}{4} r^2 + 2r - 2 = 0.
$$
Its positive root $r_0 = 0.6966...$.

\bigskip\bigskip
\centerline{\bf 3. THE SECOND PROBLEM}

\bigskip\noindent 
{\bf 1}. \ For a long time, the result of Tukia mentioned in the 
introduction remained the only known fact
connecting Teichm\"{u}ller spaces with geometric convexity. Recently
it was established in \cite {Kr6} that the universal Teichm\"{u}ller
space $\T$ cannot be mapped biholomorphically onto a bounded
convex domain in a uniformly convex Banach space, in particular,
onto a convex domain in the Hilbert space. This yields a restricted
negative answer to Sullivan's question.

The uniform convexity of a Banach space $X$ means strong convexity of its
unit ball; namely, for any $x_n, \ y_n$ satisfying
$\|x_n\| \le 1, \ \|y_n\| \le 1, \ \|x_n + y_n\| \to 2$  must be
$\|x_n - y_n\| \to 0$.
The uniformly convex spaces are reflexive and have another important
property: any bounded subset $E \subset X$ is weakly compact.
Moreover, if a sequence
$\{x_n\} \subset X$ is weakly convergent to $x_0$ and
$\|x_n\| \to \|x_0\|$ , then $x_n \to x_0$ in strong topology
of the space $X$ induced by its norm.
All this is valid, for example, for any Hilbert space and for $L_p$
spaces with $p > 1$.

\bigskip\noindent 
{\bf 2}. \ 
As for the finite dimensional case, we can show that the answer is 
{\bf negative} for the spaces $\T(0, n)$ of the punctured spheres 
(the surfaces of genus zero). Let  
$$ 
\C_{\mathbf a} = \wh \C \setminus \{a_1, \dots, a_n\}, 
\quad \hC = \C \cup \{\infty\}, 
$$
where $\mathbf a = (a_1, \dots, a_n)$ is an ordered
collection of $n > 4$ deleted distinct points. Note that $\dim \T(0,
n) = n - 3$ and that the one-dimensional space $\T(0, 4)$ is
conformally equivalent to the disk.

\bigskip\noindent
{\bf Theorem 2}. {\em There is an integer $n_0 > 4$ such that any space
$\T(0, n)$ with $n \ge n_0$ cannot be mapped biholomorphically onto
a bounded convex domain in $\C^{n-3}$.}

\bigskip
The proof of this theorem also involves conformally rigid domains 
(as for all results mentioned above) and an important 
interpolation theorem for bounded univalent functions
in the plane domains. 
This approach also has other interesting applications that are not 
presented here.

\bigskip\noindent
{\bf 3}. \ 
First we recall some needed facts from the Teichm\"{u}ller space
theory.
Consider the ordered $n$-tuples of points
\begin{equation}\label{5}
\mathbf a = (1, i, - 1, a_1, \dots, a_{n-3}), \quad n > 4,
\end{equation}
with distinct $a_j \in \hC \setminus \{1, i, -1\}$ and the corresponding
punctured spheres
$$
X_{\mathbf a} = \hC \setminus \{1, i, - 1, a_1 ,\dots, a_{n-3}\},
\quad \hC = \C \cup \{\infty\},
$$
regarded as the Riemann surfaces of genus zero. Fix a collection
$$ 
\mathbf a^0 = (1, i, - 1, a_1^0, \dots, a_{n-3}^0)
$$ 
defining the
base point $X_{\mathbf a^0}$ of Teichm\"{u}ller space $\T(0, n) =
\T(X_{\mathbf a^0})$. Its points are the equivalence classes $[\mu]$
of Beltrami coefficients from the ball
$$
\Belt(\C)_1 = \{\mu \in L_\infty(\C): \ \|\mu\|_\infty < 1\},
$$
under the relation that $\mu_1 \sim \mu_2$ if the corresponding
quasiconformal homeomorphisms $w^{\mu_1}, w^{\mu_2}: \ X_{\mathbf
a^0} \to X_{\mathbf a}$ (the solutions of the Beltrami equation
$\overline{\partial} w = \mu \partial w$ with $\mu = \mu_1, \mu_2$)
are homotopic on $X_{\mathbf a^0}$ (and hence coincide in the points
$1, i, - 1, a_1^0, \dots, a_{n-3}^0$). This models $\T(0, n)$ as the
quotient space
$\T(0, n) = \Belt(\C)_1/\sim$ 
with complex Banach structure of dimension $n - 3$ inherited from
the ball $\Belt(\C)_1$. Note that $\T(0, n)$ is a complete metric
space with intrinsic Teichm\"{u}ller metric defined by
quasiconformal maps. By Royden's theorem, this metric equals the
Kobayashi metric determined by the complex structure.

Another canonical model of $\T(0, n) = \T(X_{\mathbf a^0})$ is
obtained using the uniformization of Riemann surfaces and the
holomorphic Bers embedding of \linebreak Teichm\"{u}ller spaces. 
We now consider the disks
$$
\Delta = \{z: \ |z| < 1\}, \quad \Delta^* = \{z \in \hC: \ |z| > 1\}
$$
and the ball of Beltrami coefficients (conformal structures on $\D$)
$$
\Belt(\Delta)_1 = \{\mu \in L_\infty(\C) : \ \mu|\Delta^* = 0, \ \| \mu
\|_\infty < 1\}.
$$ 
and model the universal Teichm\"{u}ller space  $\T =\T(\D)$ as 
the space of quasi-\linebreak symmetric homeomorphisms of 
the unit circle $S^1 =
\partial \Delta$ factorized by M\"{o}bius maps. The canonical complex
Banach structure on $\T$ is defined by factorization of this ball,
letting $\mu, \nu \in \Belt(\Delta)_1$ be equivalent if the
corresponding quasiconformal maps $w^\mu, / w^\nu$ of $\hC$ coincide
on the circle $S^1$ and passing to their Schwarzian derivatives 
$S_{w^\mu}(z)$ in $D^*$ now running over a bounded
domain in the space $\B = \B(\Delta^*)$ of holomorphic 
functions $\varphi$ in $\Delta^*$ with norm
$\|\varphi\|= \sup_{\D^*} (|z|^2 - 1)^2|\varphi(z)|$. This domain is
contained in the ball $\{\|\varphi\|_\B < 1/6\}$.

The map  $\phi: \ \mu \to S_{w^\mu}$ is holomorphic and descends to
a biholomorphic map of the space $\T$ onto this domain, which we
will identify with $\T$. It contains as complex submanifolds the
Teichm\"{u}ller spaces of all hyperbolic Riemann surfaces and of
Fuchsian groups.

As is well-known, the space $\T$ coincides with the union of inner
points of the set
$$
\mathbf S =  \{\varphi = S_w \in \B: \ w \ \text{univalent in} \
\Delta^* \};
$$
on the other hand, by Thurston's theorem, $\mathbf S \setminus
\overline{\T}$ has uncountable many isolated points  
$\varphi_0 = S_{w_0}$
which correspond to conformally rigid domains $w_0(\Delta^*)$.

\bigskip\noindent
{\bf 4}. \  
Using the holomorphic universal covering map $h: \
\Delta \to X_{\mathbf a^0}$, one represents the surface $X_{\mathbf
a^0}$ as the quotient space $\Delta/\Gamma_0$ (up to conformal
equivalence), where $\Gamma_0$ is a torsion free Fuchsian group of the
first kind acting discontinuously on $\Delta \cup \Delta^*$. The
functions $\mu \in L_\infty(X_{\mathbf a^0})$ are lifted to $U$ as the
Beltrami $(-1, 1)$-measurable forms  $\wt \mu d\overline{z}/dz$ in
$\Delta$ with respect to $\Gamma_0$ which satisfy $(\wt \mu \circ \gamma)
\overline{\gamma^\prime}/\gamma^\prime = \wt \mu, \ \gamma \in \Gamma_0$ 
and form the
Banach space $L_\infty(\Delta, \Gamma_0)$.

We extend these $\wt \mu$ by zero to $U^*$ and consider the unit
ball $\Belt(\Delta, \Gamma_0)_1$ of $L_\infty(\Delta, \Gamma_0)$. Then the
corresponding Schwarzians $S_{w^{\wt \mu}|\Delta^*}$ belong to the
universal Teichm\"{u}ller space $\T$ and the subspace of such
Schwarzians is regarded as the {\em Teichm\"{u}ller space
$\T(\Gamma_0)$ of the group $\Gamma_0$}. It is canonically
isomorphic to the space $\T(X_{\mathbf a^0})$. Moreover,
 \begin{equation}\label{6}
\T(\Gamma_0) = \T \cap \B(\Gamma_0),
\end{equation}
where $\B(\Gamma_0)$ is an $(n - 3)$-dimensional subspace of $\B$
which consists of elements $\varphi \in \B$ satisfying
$(\varphi \circ \gamma) (\gamma^\prime)^2 = \varphi$ for all 
$\gamma \in \Gamma_0$; see, e.g., \cite{Le}. 
This leads to the representation of the space
$\T(X_{\mathbf a^0})$ as a bounded domain in the complex Euclidean
space $\C^{n-3}$.

Note also that the space $\B$ is dual to the subspace
$A_1(\Delta^*)$ in $L_1(\Delta^*)$ formed by integrable holomorphic
functions in $\Delta^*$, while $B(\Delta^*, \Gamma_0)$ has the same
elements as the space $A_1(\Delta^*, \Gamma_0)$ of integrable
holomorphic forms of degree $- 4$ with norm $\|\varphi\| =
\iint_{\Delta^*/\Gamma_0} |\varphi(z)| dx dy$.

\bigskip\noindent
{\bf 5}. \ The collections (5) fill a domain $\mathcal U_n$ in
$\C^{n-3}$ obtained by deleting from this space the hyperplanes $\{z
= (z_1, \dots, z_{n-3}): \ z_j = z_l, \ j \ne l\}$, and with $z_1 =
1, \ z_2 = i, \ z_3 = - 1$. This domain represents the Torelli space
of the spheres $X_{\mathbf a}$ and is covered by $\T(0, n)$. Namely,
we have (see, e.g., [14, Section 2.8])

\bigskip\noindent
{\bf Lemma 1}. {\em The holomorphic universal covering space of
$\mathcal U_n$ is the Teichm\"{u}ller space $\T(0, n)$. This means
that for each punctured sphere $X_{\mathbf a}$, there is a
holomorphic universal covering
$$
\pi_{\mathbf a}: \T(0, n) = \T(X_{\mathbf a}) \to \mathcal U_n.
$$
The covering map $\pi_a$ is well defined by
$$
\pi_{\mathbf a} \circ \phi_{\mathbf a}(\mu) =
(1, i, - 1, w^\mu(a_1), \dots, w^\mu(a_{n-3})),
$$
where $\phi_{\mathbf a}$ denotes the canonical projection of the
ball $\Belt(\Delta)_1$ onto the space $\T(X_{\mathbf a})$. }

Truncated collections $\mathbf a_{*} = (a_1, \dots, a_{n-3})$
provide the local complex coordinates on the space $\T(0, n)$ and
define its complex structure.

\bigskip
Let us consider the ball $\Belt(\Delta)_1$ and call its elements $\mu$
defining the same point of the universal Teichm\"{u}ller space {\em
$\T$-equivalent}. The corresponding homeomorphisms $w^\mu$ coincide
on the unit circle.

We now assume that the coordinates $a_j^0$ of the surface
$X_{\mathbf a^0}$ are placed on the circle $S^1$ and define on this
ball another equivalence relation, letting $\mu, \ \nu \in
\Belt(\Delta)_1$ be equivalent if $w^\mu(a_j^0) =w^\nu(a_j^0)$ for
all $j$ and the homeomorphisms $w^\mu, \ w^\nu$ are homotopic on the
punctured sphere $X_{\mathbf a^0}$. Let us call such $\mu$ and $\nu$
{\em strongly $n$-equivalent}. This equivalence is weaker 
than $\T$- equivalence, i.e., if two coefficients $\mu, \nu \in
\Belt(\Delta)_1$ are $\mathbf T$-equivalent, then they are also
strongly $n$-equivalent, which implies, (by descending to the 
equivalence classes) a holomorphic map $\chi$ of the underlying 
space $\T$ into
$\T(0, n) = \T(X_{\mathbf{a^0}})$.

This map is a split immersion, i.e., it has local holomorphic
sections. In fact, we have much more:

\bigskip\noindent
{\bf Lemma 2}. {\em The map $\chi$ is surjective and has a global
holomorphic section \linebreak  $s: \T(X_{\mathbf{a^0}}) \to \T$. }

\bigskip\noindent
\textbf{Proof}. The surjectivity of $\chi$ is a consequence of the
following interpolation result from \cite{CHMG}.

\bigskip\noindent
{\bf Lemma 3}. {\em Given two cyclically ordered collections of
points $(z_1, \dots, z_m)$ and $(\zeta_1, \dots, \zeta_m)$ on the unit
circle $S^1 = \{|z| = 1\}$, there exists a holomorphic univalent
function $f$ in the closure of the unit disk $\Delta = \{|z| < 1\}$
such that $|f(z)| < 1$ for $z \in \overline{\Delta}$ distinct from $z_1,
\dots, z_m$, and $f(z_k) = \zeta_k$ for all $k = 1, \dots, m$.
Moreover, there exist univalent polynomials $f$ with such an
interpolation property. }

\bigskip
Since the interpolating function $f$ given by this lemma is regular
up to the boundary, it can be extended quasiconformally across the
boundary circle $S^1$ to the whole sphere $\hC$. Hence, given a
cyclically ordered collection $(z_1, \dots, z_m)$ of points on
$S^1$, then for any ordered collection $(\zeta_1, \dots, \zeta_m)$ in
$\hC$, there is a quasi-\linebreak conformal homeomorphism 
$\wh f$ of the whole
sphere $\hC$ carrying the points $z_j$ to $\zeta_j, \ j = 1, \dots, m$,
and such that its restriction to the closed disk $\overline \Delta$ is
biholomorphic on $\overline{\Delta}$ (and similarly for the ordered
collections of points on arbitrary quasicircles).

Applying Lemma 1, one constructs quasiconformal extensions of $f$
lying in prescribed homotopy classes of homeomorphisms $X_{\mathbf
z} \to X_{\mathbf w}$.

\bigskip
To prove the assertion of Lemma 2 on holomorphic sections for
$\chi$, take a dense subset
$$
e = \{x_1, \ x_2, \ \dots\} \subset X_{\mathbf a^0} \cap S^1
$$
accumulating to all points of $S^1$ and consider the surfaces
$$
X_{\mathbf a^0}^m = X_{\mathbf a^0} \setminus \{x_1, \dots, x_m\},
\quad m \ge 1
$$
(having type $(0, n + m)$). The equivalence relations on
$\Belt(\C)_1$ for $X_{\mathbf a^0}^m$ and $X_{\mathbf a^0}$ generate
a holomorphic map
$\chi_m: \T(X_{\mathbf a^0}^m) \to \T(X_{\mathbf a^0})$. 

The inclusion map $j_m: \ X_{\mathbf a^0}^m \hookrightarrow
X_{\mathbf a^0}$ forgetting the additional punctures generates a
holomorphic embedding $s_m: \ \T(X_{\mathbf a^0}) \hookrightarrow
\T(X_{\mathbf a^0}^m)$ inverting $\chi_m$.
To present this section analytically, we uniformize the surface
$X_{\mathbf a^0}^m$ by a torsion free Fuchsian group $\Gamma_0^m$ on
$\Delta \cup \Delta^*$ so that $X_{\mathbf a^0}^m =
\Delta/\Gamma_0^m$. By (6), its Teichm\"{u}ller space
$\T(\Gamma_0^m) = \T \cap \B(\Gamma_0^m)$.

The holomorphic universal covering maps $h: \ \Delta^* \to
\Delta^*/\Gamma_0$ and $h^m: \ \Delta^* \to \Delta^*/\Gamma_0^m$ are
related by $j \circ h^m = h \circ \wh j$, where $\wh j$ is the lift
of $j$. This induces a surjective homomorphism of the covering
groups $\theta_m: \Gamma_0^m \to \Gamma_0$ by
 \begin{equation}\label{7}
 \wh j \circ \gamma = \theta_m(\gamma) \circ \gamma, 
\quad \gamma \in \Gamma_0^m,
\end{equation}
and the norm preserving isomorphism $\wh j_{m,{*}}: \ \B(\Gamma_0)
\to \B(\Gamma_0^m)$ by
\begin{equation}\label{8}
\wh j_{m,{*}} \varphi = (\varphi \circ \wh j) (\wh j^\prime)^2,
\end{equation}
which projects to the surfaces $X_{\mathbf a^0}$ and $X_{\mathbf
a^0}^m$ as the inclusion of the space $Q(X_{\mathbf a^0})$ of
quadratic differentials corresponding to $\B(\Gamma_0)$ into the
space $Q(X_{\mathbf a^0}^m)$ (cf. \cite{EK}). The equality (8)
represents the section $s_m$ indicated above.

\bigskip\noindent
{\bf 6}. \  
To investigate the limit function for $m \to \infty$, we embed $\T$
into the space $\B$ and compose each $s_m$ with a biholomorphism
$$
\eta_m: \ \T(X_{\mathbf a^0}^m) \to \T(\Gamma_0^m) = \T \cap
\B(\Gamma_0^m) \quad (m = 1, 2, \dots).
$$
Then the elements of $\T(\Gamma_0^m)$ are represented in the form
$$
\wh s_m(z, \cdot) = S_{f^m}(z; X_{\mathbf a}),
$$
being parameterized by the points of $\T(X_{\mathbf a^0})$.

Each $\Gamma_0^m$ is the covering group of the universal cover $h_m:
\ \Delta* \to X_{\mathbf a_0^m}$, which can be normalized
(conjugating appropriately $\Gamma_0^m$) by $h_m(\infty) = \infty, \
h_m^\prime(\infty) > 0$. Take its fundamental polygon $P_m$ obtained as
the union of the regular circular $m$-gon in $\Delta^*$ centered at
the infinite point with the zero angles at the vertices and its
reflection with respect to one of the boundary arcs. These polygons
increasingly exhaust the disk $\Delta^*$ from inside; hence, by the
Carath\'{e}odory kernel theorem, the maps $h_m$ converge to the
identity map locally uniformly in $\Delta^*$.

Since the set of punctures $e$ is dense, it completely determines
the equivalence classes $[w^\mu]$ and $S_{w^\mu}$ of $\T$, and the
limit function $s(z, \cdot) = \lim_{m\to \infty} \wh s_m(z, \cdot)$
maps $\T(X_{\mathbf a^0})$ into $\T$. For any fixed $X_{\mathbf a}$,
this function is holomorphic on $\Delta^*$; hence, by the well-known
property of elements in the functional spaces with sup-norms, $s(z,
\cdot)$ is holomorphic also in the norm of $\B$. This $s$ determines
a holomorphic section of the original map $\chi$, which completes
the proof of Lemma 3. 

\bigskip\noindent 
{\bf 7}. \ 
The following lemma is a special case of the general approximation
lemma in \cite {Kr3}; it reveals some special features which are 
used also in the proof of Theorem~ 2. 

{\bf Lemma 4}. {\em For any Schwarzian $\varphi \in \T$ holomorphic in
the disk $\Delta_r^* = \{|z| > r\}, \ r < 1$, there exist a sequence
of torsion free Fuchsian groups $\Gamma_m^r$ of the first kind
acting on $\Delta_r^*$, which does not depend on $\varphi$, and a
sequence of elements $\varphi_m \in \T(\Gamma_m^r)$ canonically
determined by $\varphi$ and converging to $\varphi$ uniformly on
$\overline{\Delta^*}$; hence, 
$\lim\limits_{m\to \infty} \|\varphi_m - \varphi\|_\B = 0$.}

\noindent
\textbf{Proof}. We pass to maps $w^\mu$ preserving the
points $0, 1, \infty$ (which  does not affect their Schwarzians
$S_{w^\mu}$ forming the space $\T$) and pick on the unit circle 
$S^1$ a dense subset of dyadic points
$$
a_l^{(n)} = e^{\pi l i/2^m}; \ l = 0, 1, \dots, 2^{m+1} - 1; \
m = 2, 3, \dots \ . 
$$
Regarding the collections
$$
\mathbf a^0(r,m) = \{0, r, \ re^{\pi l i/2^{m-3}}, \infty; \ 
l = 0, 1, \dots, m - 1\}
$$
as the punctures of the base points $X_{\mathbf a^0(r,n)}$
of the spaces  $\T(0, m) = \T(X_{\mathbf a^0(r,m)})$, 
consider for each $m$ the covering group $\Gamma_m^r$ of the
universal cover $h_m: \ \Delta_r^* \to X_{\mathbf a_0^{(r,m)}}$ with
$h_m(\infty) = \infty, \ h_m^\prime(\infty) > 0$ and take its canonical
fundamental polygon $P_m$ in $\Delta_r^*$ centered at the infinite
point with the zero angles at the vertices. These polygons
increasingly exhaust the disk $\Delta_r^*$ from inside, hence the
maps $h_m$ converge to the identity map locally uniformly in
$\Delta_r^*$.

The classical result of geometric function theory implies that for
each non-zero $\varphi \in \B(\Delta_r^*)$ and large 
$m \ge m_0(\varphi)$,
the corresponding $\Gamma_m^r$-quadratic differentials
\begin{equation}\label{9} 
\varphi_m(z) = \sum\limits_{\gamma\in \Gamma_m^r}
\varphi(\gamma z) \gamma^\prime (z)^2
\end{equation}
also do not vanish and are the Schwarzians of univalent functions
$w_m$ on $\Delta_r^*$ compatible with these groups. 
The sequences $\{\Gamma_m^r\}$ and $\{\varphi_n\}$ satisfy the assertion
of the lemma.

\bigskip
Now, to complete the proof of Theorem 2, assume, to the contrary, 
that there exists an infinite sequence
of spaces $\T(0, n)$ admitting biholomorphic homeomorphisms $\eta_n$
onto the bounded convex domains $D_n \subset \C^{n-3}$, where $n$
runs over an infinite subsequence from $\field{N}$.

We embed these domains $D_n$ biholomorphically as convex
submanifolds $V_n$ into the unit ball $B(l^2)$ of the Hilbert space
$l^2$ of sequences so that each $V_n$ is
placed in $n - 3$-dimensional subspace $l_n^2$ of $l^2$ formed by
points $\mathbf c = (c_j)$ with $c_j = 0$ for all $j > n$ and
contains its origin,
$V_n \subset V_{n+1}$, and $V_n$ touches $V_{n+1}$ from inside
in its boundary point $\mathbf c_n \in \partial V_n$ whose
distance from the origin is maximal. Their union
\begin{equation}\label{10}
V_\infty = \bigcup_n V_n
\end{equation}
is a convex submanifold in the ball $B(l^2)$ whose completion $\wh
V_\infty$ is a convex domain $\wh V_\infty$ in a subspace $l_0^2$ 
of $l^2$.

Now take a Schwarzian $\varphi_{*} = S_{w^*} \in \mathbf S \setminus \T$
defining an isolated point of $\mathbf S$ (hence a conformally rigid
domain $w^*(\Delta^*)$ in $\hC$) and consider the  homotopy
functions $w_t^*(z) = t w^*(z/t)$. Each $w_t^*$ is conformal in the
wider disk $\Delta_{1/|t|}$. Pick a sequence  of positive numbers
$t_j$ approaching $1$ and apply Lemma 6 to approximate each
Schwarzian
$$
\psi_j(z) := S_{w_{t_j}^*}(z) = t_j^{-2} S_{w^*}(z/t_j)
$$
by differentials $\varphi_{m_j} \in \T(\Gamma_m^{t_j})$ satisfying
$$
\|\varphi_{m_j} - \psi_j||_\B < \frac{1}{2^j} \dist (\psi_j, \partial \T).
$$
These $\varphi_{m_j}$ are determined by $\psi_j$ (hence by original
$\varphi_{*}$) via (9) and are convergent to $\varphi_{*}$ 
locally uniformly in $\Delta^*$.

Moreover, the proof of Lemma 6 shows that one can choose in the
series (9) a sufficiently large number $m_j = n$ so that
$\T(\Gamma_m^{t_j}) = \T(\Gamma_n)$ is one of the spaces listed
above equivalent to convex domains $V_n$.

We have for each $n$ commutative diagram
$$
\CD
\T(\Gamma_{n+1})   @>\chi_{n,n+1}>>     \T(\Gamma_n)  \\
@V{\eta_{n+1}}VV                  @VV{\eta_n}V   \\
V_{n+1}         @>\wt \chi_{n,n+1}>>         V_n
\endCD
$$
where $\chi_{n,n+1}$ is again a holomorphic map generated by
forgetting the additional puncture on the base point of $\T(0, n)$
and $\wt \chi_{n,n+1} = \eta_n^{-1} \circ \chi_{n,n+1} \circ
\eta_{n+1}$. We can replace in (10) each domain $V_n$ by its image
$\wt \chi_{n,n+1}^{-1}(V_n)$ in $V_{n+1}$.

Denote by $\chi_n$ the holomorphic map $\T \to \T_n$ given by Lemma
3. Its composition with $\eta_n$ tends as $n \to \infty$ to a
holomorphic map
$$
\eta_\infty = \lim\limits_{n\to \infty} \eta_n \circ \chi_n: \ 
\T \to V_\infty.
$$
Its holomorphy is ensured by the infinite dimensional analog of
Montel's theorem following from the Alaoglu-Bourbaki theorem.
 
It follows from Lemma 3 that $\eta_\infty$ has a holomorphic
section $\sigma_\infty: \ V_\infty \to \T$ mapping $V_\infty$
biholomorphically onto a domain
$$
\T_\infty = \bigcup_n \chi_{n,n+1}^{-1} \T(\Gamma_n) \subset \T \cap
\B_0,
$$
where $\B_0$ is some subspace of $\B$ (cf. (6)). Its inverse 
$\eta_\infty^{-1}$ also is holomorphic. 

Noting that the sequence of images 
$\x_n = \eta_\infty(\varphi_n) \in V_\infty$ is weakly
compact in the space $l^2$ and passing if needed to a convergent
subsequence to some point $\x_0 \in l^2$, one gets
 \begin{equation}\label{11}
\|\x_0\|_{l^2} \le \lim\limits_{n\to \infty} \|\x_n\|_{l^2}.
\end{equation}
Our goal is to show that only the equality is possible here, i.e.,
$\|\x_0\|_{l^2} = \lim\limits_{n\to \infty} \|\x_n\|_{l^2}$. To this
end, we consider the space $l_0^2$ as a real space with the same
norm (admitting multiplication of $\x \in l_0^2$ only with $c \in
\R$). Denote this real space by $\wt l_0^2$. The domain $V_\infty$ is
convex in $\wt l_0^2$; thus its Minkowski functional
$$
\alpha(\x) = \inf\{t > 0: \ t^{-1} \x \in V_\infty\} \quad (\x \in \wt
l_0^2)
$$
determines on this space a norm equivalent to initial norm
$\|\x\|_{l^2}$. Denote the space with the new norm by $\wt l^2_\alpha$
and notice that the domain $V_\infty$ is its unit ball.

The sequence ${x_n}$ is weakly convergent also on $\wt l^2_\alpha$;
thus, similar to (11),
$$
\alpha(\x_0) \le \lim\limits_{n\to \infty} \alpha(\x_n) \le 1.
$$
This implies that the point $\x_0$ belongs to the closure of the domain
$V_\infty$ in $l^2$-norm.

If $\alpha (\x_0) < \lim\limits_{n\to \infty} \alpha(\x_n)$ or
$\alpha(\x_0) = \lim\limits_{n\to \infty} \alpha(\x_n) < 1$, in both
these cases the point $\x_0$ must lie inside $V_\infty$. Then its
inverse image $\eta^{-1}(\x_0) \in \T$ and thus is the Schwarzian
$S_{w_0}$ of some univalent function $w_0$ on $\Delta^*$. Since
$\eta^{-1}(\x_n) = \varphi_n$ are convergent locally uniformly on
$\Delta^*$ to $S_{w^*}$, it must be $w_0 = w^*$ which yields that
$S_{w^*}$ must lie in $\T_\infty \subset \T$, in contradiction to that
it is an isolated point of the set $\mathbf S$.

It remains the case $\alpha(\x_0) = \lim\limits_{n\to \infty}
\alpha(\x_n) = 1$ which is equivalent to
 \begin{equation}\label{12}
\lim\limits_{n\to \infty} \|\mathbf{x}_n\|_{l^2} =
\|\mathbf{x}_0\|_{l^2} \quad \text{and} \ \ \mathbf{x}_0 \in
\partial V_\infty.
\end{equation}
The weak convergence $\mathbf{x}_n \to \mathbf{x}_0$ in $l^2$ and
the equality (12) together imply the strong convergence
$\lim\limits_{n\to \infty} \|\mathbf{x}_n - \mathbf{x}_0\|_{l^2} = 0$.

Then, since $\eta_\infty$ is a biholomorphic homeomorphism, the inverse
images $\eta_\infty^{-1}(\mathbf{x}_n) = \varphi_n$ must approach the
boundary of $\T$ in $\B$ and therefore $S_{w*}$ must be a boundary
point of $\T$, again contradicting that it is an isolated point of
$\mathbf S$. This completes the proof of the theorem.

\medskip
{\small\em{ \leftline{Department of Mathematics, Bar-Ilan
University,} \leftline{52900 Ramat-Gan, Israel} \leftline{and
Department of Mathematics, University of Virginia,}
\leftline{Charlottesville, VA 22904-4137, USA}}}

\end{document}